\documentclass[oneside,english]{amsart}
\usepackage[T1]{fontenc}
\usepackage[latin9]{inputenc}
\usepackage{pifont}
\usepackage{amsthm}
\usepackage{amssymb}

\makeatletter

\newcommand{\noun}[1]{\textsc{#1}}

\numberwithin{equation}{section}
\numberwithin{figure}{section}
  \theoremstyle{plain}
  \newtheorem*{thm*}{\protect\theoremname}
\theoremstyle{plain}
\newtheorem{thm}{\protect\theoremname}
  \theoremstyle{plain}
  \newtheorem{cor}[thm]{\protect\corollaryname}

\makeatother

\usepackage{babel}
  \providecommand{\corollaryname}{Corollary}
  \providecommand{\theoremname}{Theorem}
\providecommand{\theoremname}{Theorem}

\begin{document}

\title{Automorphism groups of Koras-Russell threefolds of the second kind}

\author{Charlie Petitjean}

\address{Charlie Petitjean, Institut de Math\'ematiques de Bourgogne, Universit\'e
de Bourgogne Franche-Comt\'e, 9 Avenue Alain Savary, BP 47870, 21078
Dijon Cedex, France}

\email{charlie.petitjean@u-bourgogne.fr}


\begin{abstract}
We determine the automorphism groups of Koras-Russell threefolds of
the second kind. In particular we show that these groups are semi-direct
products of two subgroups, one given by the multiplicative group and
the other isomorphic to a polynomial ring in two variables with the
addition law. We also show that these groups are generated by algebraic
subgroups isomorphic to $\mathbb{G}_{m}$ and $\mathbb{G}_{a}$.
\end{abstract}

\maketitle

\section{Introduction}

In this article, we study automorphism groups of a family of affine
threefolds, called \emph{Koras-Russell threefolds of the second kind}.
These threefolds are smooth, contractible affine varieties which first
appeared in the work of Koras and Russell on the linearization problem
for action of the multiplicative group $\mathbb{G}_{m}$ on $\mathbb{A}^{3}$(
see \cite{Ka-K-ML-R,Ka-ML,K-R,ML1}. They can be roughly classified
into three types according to the richness of the additive group action
$\mathbb{G}_{a}$ on them.

Koras-Russell threefolds of the first kind are defined by equations
of the form: $\{x+x^{d}y+z^{\alpha_{2}}+t^{\alpha_{3}}=0\}$ in $\mathbb{A}^{4}=\mathrm{Spec}(\mathbb{C}[x,y,z,t])$,
where $2\leq d$, $2\leq\alpha_{2}\leq\alpha_{3}$ with $\mathrm{gcd}(\alpha_{2},\alpha_{3})=1$.
These varieties admit several $\mathbb{G}_{a}$-actions. 

Koras-Russell threefolds of the second kind are defined by equations
of the form: $\{x+y(x^{d}+z^{\alpha_{2}})^{l}+t^{\alpha_{3}}=0\}$
in $\mathbb{A}^{4}=\mathrm{Spec}(\mathbb{C}[x,y,z,t])$, where $2\leq d$,
$1\leq l$, $2\leq\alpha_{2}\leq\alpha_{3}$ with $\mathrm{gcd}(\alpha_{2},d)=\mathrm{gcd}(\alpha_{2},\alpha_{3})=1$.These
varieties admit essentially one $\mathbb{G}_{a}$-action ( see the
second section). Such varieties are called \emph{semi-rigid} (see
\cite[4.1.2]{A}).

The other Koras-Russell threefolds which will be called of third kind
do not admit non trivial actions of the additive group; such varieties
are called rigid.

Although initially used only as part of the proof of linearization
of $\mathbb{G}_{m}$-actions on $\mathbb{A}^{3}$, these varieties
are now of interest in their own right as exotic hypersurfaces and
also as counter-example of the cancellation property for example (
see \cite{D,D-MJ-P,MJ,D-MJ-P1}).

By design all Koras-Russell threefolds admit hyperbolic $\mathbb{G}_{m}$-actions
with a unique fixed point. The questions considered here are: what
is the structure of their automorphism groups, what are the actions
of the additive and multiplicative groups and more generally what
are the algebraic subgroups.

The problem of describing the group of polynomial automorphisms of
affine spaces is a classical subject in algebraic geometry. In dimension
two, the theorem of Jung and Van der Kulk \cite{J,VK} shows that
every automorphism of the polynomial ring in two variables can be
decomposed into the product of affine automorphisms and de Jonquieres
automorphisms and in addition we have a structure of amalgamated product
of these two groups. However, there is no similar theorem to describe
the automorphism groups of the polynomial ring in $n$ variables for
$n\geq3$. It is therefore interesting to have a description of the
automorphism groups of three dimensional algebraic varieties which
resemble $\mathbb{A}^{3}$, in the sense that they are smooth rational
and contractible.

We recall that the determination of the automorphism groups of an
affine variety $X$ is equivalent to that of its ring of regular functions
$\mathbb{C}[X]$. We denote for simplicity by $Aut(X)$ the group
$Aut_{\mathbb{C}}(\mathbb{C}[X])$ of $\mathbb{C}$-automorphisms
of $\mathbb{C}[X]$.

The natural correspondence between an affine variety and its ring
of regular functions gives a particular interpretation for $\mathbb{G}_{a}$-actions.
Indeed, the set of all $\mathbb{G}_{a}$-actions on $X$ is in one-to-one
correspondence with the set of all \emph{locally nilpotent derivations}
on $\mathbb{C}[X]$ (see \cite{F,ML2}). We recall that a derivation
$\partial$ on $\mathbb{C}[X]$ is called a locally nilpotent derivation\emph{
}if for any $f\in\mathbb{C}[X]$, there exists $n\in\mathbb{Z}_{\geq0}$
such that $\partial^{n}(f)=0$. The set of all locally nilpotent derivations
on $\mathbb{C}[X]$ is denoted by $\mathrm{LND}(\mathbb{C}[X])$\emph{.}
The intersection of the kernels of all locally nilpotent derivation
on $\mathbb{C}[X]$, $\mathrm{ML}(X):=\bigcap_{\partial\in\mathrm{LND}(\mathbb{C}[X])}ker(\partial)$,
is called the \emph{Makar-Limanov invariant} of $X$. It is a subring
of the ring of regular functions $\mathbb{C}[X]$ which is invariant
by all automorphisms. When this subring is non trivial, that is, not
equal to $\mathbb{C}$ or $\mathbb{C}[X]$, this fact can be used
to study the automorphism groups. 

In the case of the affine space it is trivial, that is, it is equal
to $\mathbb{C}$, thus the Makar-Limanov invariant does not provide
any information on the automorphism group of $\mathbb{A}^{n}$. The
Koras-Russell threefolds of the third kind admit a Makar-Limanov invariant
given by their ring of regular functions $\mathbb{C}[X]$ \cite[8.3]{Ka-ML}.
Once again the approach of using the Makar-Limanov does not give any
information on the automorphism groups.

The study of the automorphism groups of varieties of the first kind
has been described by Moser-Jauslin in \cite{MJ}. In this case, the
Makar-Limanov invariant is given by $\mathbb{C}[x]$ \cite[8.3]{Ka-ML},
thus these varieties are not isomorphic to $\mathbb{A}^{3}$ but admit
several $\mathbb{G}_{a}$-actions. We have the following sequence
of inclusion of rings:

\begin{flushleft}
\[
\mathbb{C}[x]\subset\mathbb{C}[x,z,t]\subset\mathbb{C}[X]\subset\mathbb{C}[x^{\pm1},z,t]\subset\mathbb{C}(x,z,t).
\]
 Moreover, the polynomial defining the hypersurface is homogeneous
for the following linear $\mathbb{G}_{m}$-action on $\mathbb{A}^{4}$:
\[
\lambda\cdot(x,y,z,t)\rightarrow(\lambda^{\alpha_{2}\alpha_{3}}x,\lambda^{-(d-1)\alpha_{2}\alpha_{3}}y,\lambda^{\alpha_{3}}z,\lambda^{\alpha_{2}}t),
\]
and so $\mathbb{G}_{m}$ is one subgroup of the automorphism groups
of these varieties. One of the results in \cite{MJ} is the following
:
\par\end{flushleft}

\begin{thm*}[MJ]Let $X$ be a Koras-Russell threefold of the first
kind. The automorphism group is isomorphic to $Aut(X)\simeq\mathcal{A}_{1}\rtimes\mathbb{G}_{m}$,
where $\mathcal{A}_{1}$ is the subgroup of \noun{$\mathbb{C}[x]$}-automorphisms
of $\mathbb{C}[x][z,t]$ which are congruent to the identity modulo
$(x)$, and which stabilize the ideal $I=(x^{d},z^{\alpha_{3}}+t^{\alpha_{2}}+x)$.

\end{thm*}

In particular, the subalgebras $\mathbb{C}[x,z,t]$ and $\mathbb{C}[x]$
are stable by any automorphism and any element of $\mathcal{A}_{1}$
can be lifted to element of $Aut(X)$.

Here we apply a similar approach for automorphism groups of Koras-Russell
threefolds of the second kind. In this case, the Makar-Limanov invariant
is equal to $\mathbb{C}[x,z]$ \cite[8.3]{Ka-ML}. In particular these
varieties are \emph{semi-rigid} (see \cite[4.1.2]{A}), that is, $\mathrm{LND}(X)=ker(\partial)\cdot\partial$
for some non trivial locally nilpotent derivation $\partial\in\mathrm{LND}(X)$.
In order to study $Aut(X)$, we first note that $X$ is rational,
and furthermore that $\mathbb{C}[X]$ is generated by four elements
$x$, $y$, $z$, $t$ which satisfy the relation $x+y(x^{d}+z^{\alpha_{2}})^{l}+t^{\alpha_{3}}=0$.
Let $f$ be the polynomial $x^{d}+z^{\alpha_{2}}$. The fraction field
of $\mathbb{C}[X]$ is generated by $\mathbb{C}(x,z,t)$ where $y=-\frac{x+t^{\alpha_{3}}}{(f)^{l}}$.
In other words any element of $Aut(X)$ is determined by the image
of $x$, $z$ and $t$. We have the following sequence of inclusion
:

\[
\mathbb{C}[x,z]\subset\mathbb{C}[x,z,t]\subset\mathbb{C}[X]\subset\mathbb{C}[x,z,t,f^{-1}]\subset\mathbb{C}(x,z,t).
\]
Moreover, the polynomial defining the hypersurface is homogeneous
for the following linear $\mathbb{G}_{m}$-action on $\mathbb{A}^{4}$:
\[
\lambda\cdot(x,y,z,t)\rightarrow(\lambda^{\alpha_{2}\alpha_{3}}x,\lambda^{-(dl-1)\alpha_{2}\alpha_{3}}y,\lambda^{d\alpha_{3}}z,\lambda^{\alpha_{2}}t),
\]
and so $\mathbb{G}_{m}$ is isomorphic to a subgroup of the automorphism
groups of these varieties. The goal of this article is to determine
completely the automorphism groups.
\begin{thm*}
Let $X$ be a Koras-Russell threefold of the second kind. The automorphism
group is isomorphic to $Aut(X)\simeq\mathcal{A}\rtimes\mathbb{G}_{m}$,
where $\mathcal{A}$ is the subgroup of $Aut(\mathbb{C}[x,z,t])$
whose elements fix $\mathbb{C}[x,z]$ and send $t$ to $t+(x^{d}+z^{\alpha_{2}})^{l}p(x,z)$
for some polynomial $p(x,z)\in\mathbb{C}[x,z]$. The group $\mathcal{A}$
is isomorphic to $(\mathbb{C}[x,z],+)$.
\end{thm*}
In particular the subalgebras $\mathbb{C}[x,z,t]$, $\mathbb{C}[x]$
(that is $\mathrm{ML}(X)$) and $\mathbb{C}[z]$ are stable by every
automorphism and every element of $\mathcal{A}$ can be lifted to
a unique element of $Aut(X)$.

\section{Proof of the Theorem}

\begin{flushleft}
Let $X$ be a Koras-Russell threefolds of the second kind given by:
\[
\{x+y(x^{d}+z^{\alpha_{2}})^{l}+t^{\alpha_{3}}=0\}\subset\mathbb{A}^{4}=\mathrm{Spec}(\mathbb{C}[x,y,z,t]).
\]
 
\par\end{flushleft}

The Makar-Limanov invariant of $X$ is equal to $\mathbb{C}[x,z]$.
We will use this fact to completely determine the actions of the additive
group $\mathbb{G}_{a}$ on $X$. Let $\partial\in\mathrm{LND}(\mathbb{C}[X])$
be the irreducible derivation given by: 
\[
\partial:=\alpha_{3}t^{\alpha_{3}-1}\frac{\partial}{\partial y}-(x^{d}+z^{\alpha_{2}})^{l}\frac{\partial}{\partial t}.
\]

\begin{flushleft}
This derivation is obtained by considering the Jacobian determinant:
$Jac_{\mathbb{C}[x,z]}(P,\cdot)=\left|\begin{array}{cc}
\frac{\partial P}{\partial y} & \frac{\partial\cdot}{\partial y}\\
\frac{\partial P}{\partial t} & \frac{\partial\cdot}{\partial t}
\end{array}\right|$with $P(x,y,z,t)=x+y(x^{d}+z^{\alpha_{2}})^{l}+t^{\alpha_{3}}$.
\par\end{flushleft}

Thus $\partial(x)=0$, $\partial(z)=0$, $\partial(t)=-(x^{d}+z^{\alpha_{2}})^{l}$
and $\partial(y)=\alpha_{3}t^{\alpha_{3}-1}$ and the kernel of $\partial$
is $\mathbb{C}[x,z]=\mathrm{ML}(X)$, preserved by the automorphisms
of $\mathbb{C}[X]$. As these hypersurfaces are semi-rigid (see \cite[4.1.2]{A})
any other locally nilpotent derivation can be written in the following
form: $\partial_{q}=q(x,z)\partial$ with $q(x,z)\in\mathrm{ML}(X)=\mathbb{C}[x,z]$.
In this case $\partial$ is the unique irreducible element of $\mathrm{LND}(X)$
up to multiplication by a constant. For any $\varphi\in Aut(X)$ ,
$\partial_{\varphi}:=\varphi^{-1}\circ\partial\circ\varphi$ is also
an irreducible element of $\mathrm{LND}(\mathbb{C}[X]$ \cite[corollary 2.3]{F}.
Therefore, there exists $a\in\mathbb{C}^{*}$ such that $\partial_{\varphi}=a\partial$
. Moreover as the kernel of $\partial$ is equal to $\mathbb{C}[x,z]$
we obtain that $\varphi(\mathbb{C}[x,z])=\mathbb{C}[x,z]$.\\
\\
\\

\ding{182} Any automorphism $\varphi$ of $X$ induces an automorphism
$\hat{\varphi}$ of $\mathbb{C}[x,z]$ satisfying:

i) the ideal $(f)$ generated by $f$ is preserved

ii) there is $\mu_{\hat{\varphi}}\in\mathbb{C}^{*}$ such that $\hat{\varphi}(x)=(\mu_{\hat{\varphi}})^{\alpha_{2}}x$
and $\hat{\varphi}(z)=(\mu_{\hat{\varphi}})^{d}z$.

The part i) comes from the fact that $\varphi$ must preserve the
Makar-Limanov invariant $ML(X)=\mathbb{C}[x,z]$, which is also the
kernel of $\partial$. Moreover, let $\pi$ be the projection on the
coordinates $(x,z)$:
\begin{equation}
\begin{array}{cccc}
\pi: & X & \rightarrow & \mathbb{A}^{2}=\mathrm{Spec}(\mathbb{C}[x_{0},z_{0}])\\
 & (x,y,z,t) & \rightarrow & (x,z),
\end{array}\label{eq:}
\end{equation}

\begin{flushleft}
and let $f_{0}:=(x_{0}^{d}+z_{0}^{\alpha_{2}})$ for any point $(x_{0},z_{0})\in\mathbb{A}^{2}$.
Then $\mbox{\ensuremath{\pi}}^{-1}(x_{0},z_{0})$ is isomorphic to
$\mathbb{A}^{1}$ if $f_{0}\neq0$ or $f_{0}=0$ and $x=0$ but it
is isomorphic to $\alpha_{3}$ copies of $\mathbb{A}^{1}$ otherwise.
Thus the cuspidal curve must be preserved by $\hat{\varphi}$, that
is, there exists $\lambda_{\hat{\varphi}}\in\mathbb{C}^{*}$ such
that $\hat{\varphi}(f)=\lambda_{\hat{\varphi}}f$.
\par\end{flushleft}

The only automorphism $\psi$ of $\mathbb{C}[x,z]$ which is congruent
to the identity modulo$f$ is indeed the identity. To see this, note
that if such an automorphism exists, it would stabilise all plane
curves defined by the level sets of $f$. Suppose that $C$ is the
zero set of the polynomial $f-c$, where $c$ is a non-zero constant.
Then $C$ is the open set of a smooth compact Riemann surface of genus
$\ge1$. By Hurwitz's theorem (see \cite{H}), the automorphism group
of $C$ is finite, and in fact the automorphism $\psi$ must be of
finite order. This implies in particular that $\psi$ is conjugate
to a linear action by \cite{J,VK}. By considering the linear part,
since $\psi$ is congruent to the identity modulo $f$, we have that
$\psi$ is conjugate to the identity, and thus is in fact the identity.

As $\hat{\varphi}\in Aut(\mathbb{C}[x,z])$ satisfies $\hat{\varphi}(f)=\lambda_{\hat{\varphi}}f$
there is $\mu_{\hat{\varphi}}\in\mathbb{C}^{*}$ such that $\hat{\varphi}(x)=(\mu_{\hat{\varphi}})^{\alpha_{2}}x$
modulo $f$ and $\hat{\varphi}(z)=(\mu_{\hat{\varphi}})^{d}z$ modulo
$f$. By composition with the linear automorphism $\varphi_{0}$:
$\varphi_{0}(x)=(\mu_{\hat{\varphi}})^{-\alpha_{2}}x$ and $\varphi_{0}(z)=(\mu_{\hat{\varphi}})^{-d}z$,
we obtain $\varphi_{0}\circ\hat{\varphi}=\psi\equiv\mathrm{id}$.
This proves part ii).\\

\ding{183} Now consider automorphisms of $\mathbb{C}[X]$ which fixed
$\mathbb{C}[x,z]$ and we focus on the image of $t$ by any element
$\varphi\in Aut(X)$. We prove first that $\varphi(t)$ is of the
form $at+h(x,z)$ with $a\in\mathbb{C}^{*}$.

As $\varphi(f)=f$ where $f=(x^{d}+z^{\alpha_{2}})$ we apply $\partial_{\varphi}$
to the variable $t$: $\partial_{\varphi}(t)=\varphi^{-1}\circ\partial\circ\varphi(t)$
thus:

\[
\begin{array}{ccc}
\varphi\circ(a\partial(t)) & = & \partial\circ\varphi(t)\\
\varphi(-af^{l}) & = & \partial\circ\varphi(t)\\
-af^{l} & = & \partial\circ\varphi(t).
\end{array}
\]

\begin{flushleft}
This means that $\partial\circ\varphi(t)=\partial(at)$ and thus $\varphi(t)-at\in ker(\partial)=\mathbb{C}[x,z]$.
\par\end{flushleft}

Secondly we prove that $a$ is a $\alpha_{3}$-th root of the unity
and $h(x,z)\in(f^{l})$. Let $J\subset\mathbb{C}[X]$ be the ideal
generated by $(f^{l})$, and let $I=J\cap\mathbb{C}[x,z,t]$. Then
$I=(f^{l},x+t^{\alpha_{3}})\subset\mathbb{C}[x,z,t]$. Indeed, $f^{l}$
and $x+t^{\alpha_{3}}$ are in $I$ thus $(f^{l},x+t^{\alpha_{3}})\subset I$.
Now suppose $Q\in I$. Since $Q\in J$, there exists $P\in\mathbb{C}[X]$
such that $Q=f^{l}P$. Any $P\in\mathbb{C}[X]$ can be decomposed
in a unique way as follows $P=\sum_{i=0}^{n}p_{i}(x,z,t)y^{i}$ such
that $f^{l}$ does not divide $p_{i}$ if $i\geq1$ since $y=-(x+t^{\alpha_{3}})f^{-l}$.
This gives $Q=f^{l}p_{0}+\sum_{i=1}^{n}p_{i}(x,z,t)(x+t^{\alpha_{3}})y^{i-1}$.
Now since $Q\in\mathbb{C}[x,z,t]$, the polynomial $p_{i}(x,z,t)=0$
for $i\geq2$ and then $Q=f^{l}p_{0}+p_{1}(x,z,t)(x+t^{\alpha_{3}})$.

Consider the ideal $I=(f^{l},x+t^{\alpha_{3}})\subset\mathbb{C}[x,z,t]$
which is preserved by $\varphi$, thus $\varphi(x+t^{\alpha_{3}})\in I$.
There exists $b(x,z,t)=\sum_{i=0}^{n}b_{i}(x,z)t^{i}$ and $c(x,z,t)=\sum_{i=0}^{m}c_{i}(x,z)t^{i}$
such that:

\begin{flushleft}
\[
\varphi(x+t^{\alpha_{3}})=x+(at+h(x,z))^{\alpha_{3}}=(\sum_{i=0}^{n}c_{i}(x,z)t^{i})(x+t^{\alpha_{3}})+(\sum_{i=0}^{m}b_{i}(x,z)t^{i})f^{l}.
\]
In addition, one can assume that for all $i$, $f^{l}$ does not divide
$c_{i}(x,z)$. Considering the highest degree in the variable $t$
on both side, that is the coefficient of $t^{\alpha_{3}}$. If $n\ge1$,
then $f^{l}$ divides $c_{n}(x,z)$, which contradicts the assumption.
Thus $n=0$, and $c_{0}(x,z)$ is congruent to $a^{\alpha_{3}}$ modulo
$f^{l}$. We can therefore suppose $c_{0}(x,z)$ equals $a^{\alpha_{3}}$,
by adding the appropriate term to $b(x,z,t)$. Moreover, the equality
above implies that $h(x,z,t)$ is in the ideal $(x,z,t)$. Now considering
the coefficient of $x$ on both sides, we see that $c_{0}(x,z)=1$.
This implies: 
\[
\varphi(x+t^{\alpha_{3}})-(x+t^{\alpha_{3}})=h(x,z)[\sum_{k=1}^{\alpha_{3}-1}\left(\begin{array}{c}
\alpha_{3}\\
k
\end{array}\right)(at)^{k}h(x,z)^{\alpha_{3}-1-k}]=\sum_{j=1}^{m}b_{j}(x,z)t^{j}f^{l},
\]
 and considering the coefficient of $t$ on both sides, we see that
$f^{l}$ divides $h(x,z)$. So $h(x,z)$ is in the ideal generated
by $f^{l}$ in $\mathbb{C}[x,z]$, that is $h(x,z)=(x^{d}+z^{\alpha_{2}})^{l}p(x,z)$.
\par\end{flushleft}

\ding{184} If $\psi\in Aut(\mathbb{C}[x,z,t])$ is in $\mathcal{A}$
then $\psi$ extends in a unique way to an automorphism in $Aut(X)$,
via the computation:

\begin{equation}
\psi(y)=\psi\left(\frac{-(x+t^{\alpha_{3}})}{f^{l}}\right)=\frac{-(x+(t+f^{l}p(x,z)){}^{\alpha_{3}})}{f^{l}}=y+H(x,z,t).\label{eq:-1}
\end{equation}
 In particular $\mathcal{A}$ can be identified with a normal subgroup
of $Aut(X)$. By part 2 and the $\mathbb{G}_{m}$-action given in
the introduction every automorphism $\varphi$ can be expressed in
a unique way as a composition of an element of $\mathcal{A}$ and
an element of the $\mathbb{G}_{m}$-action. In conclusion $Aut(X)=\mathcal{A}\rtimes\mathbb{G}_{m}$. 

\noindent\\
\\
\noindent

We recall several results of \cite[proposition 3.1]{D-MJ-P1} in order
to compare with that obtained as corollaries of the theorem of this
article. In the case of Koras-Russell threefolds of the first kind
the following results hold: every automorphism of $X$ extends to
an automorphism of $\mathbb{A}^{4}$, the group $Aut(X)$ admits 4
orbits with in particular the origin $(0,0,0,0)$ as fixed point and
the subgroup generated by all $\mathbb{G}_{a}$ and $\mathbb{G}_{m}$
actions is strictly smaller than $Aut(X)$.

In the case of Koras-Russell threefolds of the second kind, let $\pi$
be the projection defined in \ref{eq:} and consider a partition of
$\mathbb{A}^{2}=\mathrm{Spec}(\mathbb{C}[x_{0},z_{0}])$ by the origin
and the $\mathbb{G}_{m}$-stable curves $C_{\alpha,\beta}=\{\alpha x_{0}^{d}+\beta z_{0}^{\alpha_{2}}=0\}\setminus(0,0)$,
with $[\alpha:\beta]\in\mathbb{P}^{1}$. Then $\pi^{-1}(0,0)$ admits
two orbits given by $(0,0,0,0)$ and the line $\{x=z=t=0\}$ minus
the point $(0,0,0,0)$. For every curve $C_{\alpha,\beta}$ if $[\alpha:\beta]\neq[1,1]$
then $\pi^{-1}(C_{\alpha,\beta})\simeq\mathbb{A}^{1}\times\mathbb{A}_{*}^{1}$
considering $y=-\frac{x+t^{\alpha_{3}}}{(f)^{l}}$ and if $[\alpha:\beta]=[1,1]$
then $x=t^{\alpha_{3}}$ thus any choice of $t$ determines $x$ and
gives $\{t^{d\alpha_{3}}+z^{\alpha_{2}}=0\}\setminus(0,0)\simeq\mathbb{A}_{*}^{1}$
and $\mathcal{A}$ acts trivially on $\pi^{-1}(C_{1,1})\simeq\mathbb{A}^{1}\times\mathbb{A}_{*}^{1}$.
\begin{cor}
The automorphism group acts on $X$ with an infinite number of orbits:

a) One fixed point $(0,0,0,0)$.

b) The line $\{x=z=t=0\}$ minus the point $(0,0,0,0)$ isomorphic
to $\mathbb{A}_{*}^{1}$.

c) An infinite number of orbits $\pi^{-1}(C_{\alpha,\beta})\simeq\mathbb{A}^{1}\times\mathbb{A}_{*}^{1}$
for $[\alpha:\beta]\neq[1,1]$.

d) An infinite number of orbits isomorphic to $\mathbb{A}_{*}^{1}$,
whose union is $\pi^{-1}(C_{1,1})$.
\end{cor}
We have proved that there was not only the Makar-Limanov invariant
$\mathbb{C}[x,z]$ which was preserved by any automorphism of $\mathbb{C}[X]$
but also $\mathbb{C}[x,z,t]$ since for any automorphism $\varphi$
of $\mathbb{C}[X]$, $\varphi(t)=at+(x^{d}+z^{\alpha_{2}})^{l}p(x,z)$
and as the image of $y$ is determined by that of the other variable.
Using the same argument as in the computation of the equation \ref{eq:-1},
we see immediately that every automorphism of $X$ comes from the
restriction of an automorphism of $Aut(\mathbb{A}^{4})$.

\begin{flushleft}
For every $\varphi\in\mathcal{A}$ we have: 
\[
\varphi(\cdot)=exp(p\partial)=\sum_{k=0}^{\infty}\frac{p(x,z)^{k}\partial^{k}(\cdot)}{k!}.
\]
In particular $\varphi$ belongs to an algebraic subgroup of $Aut(X)$
isomorphic to $\mathbb{G}_{a}$. Thus the following holds.
\par\end{flushleft}
\begin{cor}
The automorphism group of $X$ is generated by $\mathbb{G}_{a}$-actions
and $\mathbb{G}_{m}$-actions.
\end{cor}

\end{document}